\documentclass{amsart}

\usepackage{amsmath, amsthm}
\usepackage{amssymb, latexsym}
\usepackage{amsfonts}

\newtheorem{theorem}[equation]{Theorem}
\newtheorem{Proposition}[equation]{Proposition}

\newtheorem{lemma}[equation]{Lemma}
\newtheorem{Lemma}[equation]{Lemma}

\newtheorem{Corollary}[equation]{Corollary}

\newtheorem{remark}[equation]{Remark}

\newcounter{com}

\newcommand{\beql}[1]{\begin{equation}\label{#1}}
\newcommand{\eeq} {\end{equation}}

\font\Aaa=msam10

\def\qed{\hbox{~~\Aaa\char'003}}

\font\Bbb=msbm10

\def\R{\hbox{\Bbb R}}

\def\Z{\hbox{\Bbb Z}}
\def\Q{\hbox{\Bbb Q}}
\def\F{\hbox{\Bbb F}}
\def\N{\hbox{\Bbb N}}

\def\A{\hbox{\Bbb A}}

\numberwithin{equation}{section}

\let\define=\def
\let\redefine=\def

\def\Aut{{\rm Aut}}

\redefine\C{{ \Bbb C }}

\def\G{\Gamma}

\redefine\D{{ \Delta }}

\define\({ \left( }
\define\){ \right) }
\define\[{ \left[ }
\define\]{ \right] }
\define\<{ \langle }
\define\>{ \rangle }
\let\ljunk=\{
\let\rjunk=\}
\redefine\{{\left\ljunk}
\redefine\}{\right\rjunk}

\redefine\.{ \diamomd }
\redefine\+{\oplus}
\redefine\.{ \cdot}
\redefine\#{\sharp}

\redefine\.{ \cdot }

        \font\Aaa=msam10

\font\Aaa=msam10

\def\qed{\hbox{~~\Aaa\char'003}}

\font\Bbb=msbm10

\def\R{\hbox{\Bbb R}}

\def\Z{\hbox{\Bbb Z}}
\def\Q{\hbox{\Bbb Q}}
\def\F{\hbox{\Bbb F}}
\def\C{\hbox{\Bbb C}}

\def\div{ \kern-.5pt\hbox{\big |} }
\def\ndiv{ {\not\kern-.5pt\hbox{\big |}\,} }
\def\ndivv{ {\not\kern+1.5pt\hbox{$\mid$}\,} }

\def\Aut{{\rm Aut}}

\def\B{^2\kern-.8pt B}
\def\G{^2\kern-.8pt G}
\def\EH{^2\kern-.8pt\hat  E}
\def\E{^2\kern-.8pt E}
\def\D{^3\kern-1pt D}
\def\FF{^2\kern-.8pt F}

\newdimen\refcodesize
\newbox\seriesbox
\def\proof{\noindent {\bf Proof.~}}

\begin{document}

\title[Invariable generation]
{Invariable generation of infinite groups}

\thanks{The authors acknowledge partial support from
ERC Advanced Grants 226135 (A.\hspace{1pt}L.) and
247034 (A.\hspace{1pt}S.,  W.\hspace{1pt}M.\hspace{1pt}K.), and ISF
grant 1117/13 (A.\hspace{1pt}L. and
A.\hspace{1pt}S.).
The first author is grateful for  the warm hospitality of the Hebrew
University while this paper was being written.}

       \author{William M. Kantor}

       \address{University of Oregon,
       Eugene, OR 97403}
       \email{kantor@uoregon.edu}

    \author{Alexander Lubotzky}
       \address{Institute of Mathematics, Hebrew University, Jerusalem 91904}
       \email{alex.lubotzky@mail.huji.ac.il}

    \author{Aner Shalev}
       \address{Institute of Mathematics, Hebrew University, Jerusalem 91904}
       \email{shalev@math.huji.ac.il}

{

}

\begin{abstract}
A subset $S$ of a group $G$ {\em invariably generates} $G$
if $G=\<s^{g(s)}\mid s\in S  \>$ for each choice of $g(s)\in G, s\in S$.
In this paper we study invariable generation of infinite groups,
with emphasis on linear groups. Our main result shows that
a finitely generated linear group is invariably generated by
some finite set of elements if and only if it is virtually solvable.
We also show that the profinite completion of an arithmetic group
having the congruence subgroup property is invariably generated by
a finite set of elements.
 \end{abstract}

\maketitle

\vspace{-8pt}

\centerline{\em Dedicated to the  memory of  \'Akos Seress}

\vspace{4pt}

\section{Introduction}

\label{Introduction}

In \cite{KLS} we studied the notion of invariable generation of finite
groups. The goal of this paper is to present some results, examples
and questions towards the study of this notion for infinite groups.

Following Dixon \cite{Di} we say that a group $G$ is
{\em invariably
generated}  by a subset $S$ of $G$ if
$G=\<s^{g(s)}\mid s\in S  \>$ for each choice of $g(s)\in G, s\in S$.
We also say that the group $G$ is IG if it is invariably generated
by some subset $S \subseteq G$, or equivalently,  if $G$ is invariably generated
by $G$; and that $G$ is FIG  if it is invariably generated
by some {\em finite} subset $S \subseteq G$.

The notion of invariable generation occurs naturally for Galois groups,
where elements are only given up to conjugacy.
IG groups were studied in a different language by Wiegold:
a group $G$ is IG if and only if it cannot be covered by a union of
conjugates of a proper subgroup, which amount to saying that in every
transitive permutation representation of $G$ on a set with more than
one element there is a fixed-point-free element.
Results on such groups can be found in \cite{W1, W2}.

In \cite{KLS} we show that a finite group $G$ is invariably generated
by at most $\log_2 |G|$ elements, and that every finite simple group is
invariably generated by two elements (the latter result is also
obtained in \cite{GM}).

We now turn to infinite groups. Which of them are FIG?
Our main result solves this problem for linear groups:

\begin{theorem}
\label{Theorem 1}
A linear group is FIG if and only if it is finitely generated
and virtually solvable.
\end{theorem}

By a well known result of Margulis and Soifer \cite{MS}, a finitely
generated linear group is virtually solvable if and only if all its
maximal subgroups are of finite index:

\begin{Corollary} A finitely generated linear group is FIG
if and only if all its maximal subgroups have finite index.
\end{Corollary}

We are not aware of a direct proof of this corollary.

The linearity assumption in Theorem 1.1 cannot be dropped:
FIG groups need not be virtually solvable.
For example, Corollary 2.7 below shows that the Grigorchuk group
(see \cite{Gr}) is FIG (in fact it is invariably generated by
its three natural generators). This follows from the fact that
the maximal subgroups of this group are all of index $2$.
The Grigorchuk group is residually finite, so we conclude that residually
finite FIG groups need not be virtually solvable.

Ol'shanski \cite{O} and Rips have constructed infinite groups $G$ in which
all proper non-trivial subgroups $H$ have order $p$ (for a given large
prime $p$). It can be arranged that these subgroups $H$ are not all
conjugate (Rips, private communication). If $H_1, H_2 \le G$ are non-conjugate
subgroups of order $p$ generated by elements $h_1, h_2$ respectively, then
$G$ is invariably generated by $h_1, h_2$, and $G$ is clearly not
virtually solvable. Unlike the previous example, this example also
shows that the linearity assumption in Corollary 1.2 is essential.

It is natural to ask which linear groups are IG. At the moment
we are unable to solve this problem.
Note that many linear groups are not IG. For example,
let $G = SL_n(\C)$. Then, using the Jordan form of matrices,
we see that every element $s \in G$ has a conjugate $s^{g(s)}$
lying in the Borel subgroup $B < G$ of upper triangular matrices.
This shows that $G$ is not IG.
A similar argument shows that, for $n>2$, the group $SL_n(\R)$
is not IG, using a parabolic subgroup of type $(2,n-2)$ instead
of a Borel subgroup.

More examples of groups which are not IG are given in Section 2
below. We also show in Proposition 2.4   that a linear
algebraic group over an algebraically closed field is IG if and only if it is
virtually solvable.

The situation over global fields is less clear. For example,
it would be nice to find out whether $SL_n(\Q)$ is an IG group.
A similar question may be asked for $SL_n(\Z)$ and for
arithmetic groups in general. In particular, is there a
correlation for such groups between being IG and having
the Congruence Subgroup Property (CSP)?

The situation is clearer for $p$-adic and adelic groups. We say that a
profinite group $G$ is invariably generated by a subset
$S \subseteq G$ if $\{s^{g(s)}\mid s\in S  \}$ generates $G$
{\em topologically} for each choice of $g(s)\in G, s\in S$.
It is easy to see that profinite groups are always IG, but
they are not necessarily FIG (even if they are finitely generated).
See Section 5 for details.

\begin{Proposition} Let $G$ be a simply connected simple Chevalley group.

(i) The adelic group $G( \widehat{\Z} )$ is FIG. In particular
the $p$-adic groups $G(\Z_p)$ are all FIG.

(ii) If $p > 3$ then the group $G(\Z_p)$ is invariably generated
by two elements.

\end{Proposition}

It is intriguing that, while arithmetic groups are not FIG
(by Theorem 1.1), their profinite completions
are often FIG. For example, let $G$ be a Chevalley group and suppose
the arithmetic group $G(\Z)$ has CSP. Then the profinite completion
$\widehat{G(\Z)}$ is isomorphic to the adelic group $G(\widehat{\Z})$,
so it is FIG by Proposition 1.3. The next result extends this
to general arithmetic groups, also in positive characteristic.

\begin{theorem}
Let $k$ be a global field of arbitrary characteristic,
$O$ its ring of integers, $T$ a finite
set of places containing all the archimedean ones. Let $G \le GL_n$
be a connected simply connected simple algebraic group defined over $k$,
and let $G(O_T): = G \cap GL_n(O_T)$.
Suppose $G(O_T)$ satisfies the Congruence Subgroup Property.
Then the profinite completion ${\widehat{G(O_T)}}$ is FIG.
\end{theorem}

In \cite{L} CSP for $G$ is shown to have various purely group-theoretic characterizations when $char(k)=0$ (e.g. $\widehat{G}$ is
boundedly generated).
There is no such known criterion when $char(k) > 0$.
Is the property ``$\widehat{G}$ is FIG'' equivalent to CSP?

To show this we need to prove that the profinite completions of
arithmetic groups without CSP are not FIG. We can show this
in some special cases, e.g. for $SL_2(\Z)$. More generally
we prove the following.

\begin{theorem} Let $G$ be any Fuchsian group. Then $\widehat{G}$
is not FIG.
\end{theorem}

The proof uses the probabilistic solution in \cite{LiSh} of
Higman's conjecture, that any Fuchsian group maps onto all
large enough alternating groups.

Some words on the structure of this paper. In Section 2 we prove
some preliminary results, and various examples are provided.
Theorem 1.1 is proved in Section 3. Section 4 is devoted to
profinite groups and contains proofs of  Proposition~1.3-Theorem~1.5.
In Section 5 we suggest some problems and directions for
further research.

We are grateful to Slava Grigorchuk, Andrei Rapinchuk and Ilya Rips
for valuable advice.

\section{Preliminary results}
\label{Preliminaries}

Let $G$ be a group and $H \le G$ a subgroup. Define
$$
\widetilde{H} = \bigcup_{g \in G} H^g,
$$
the union of all conjugates of $H$ in $G$.

The following is straightforward.

\begin{lemma}
 A subset $S \subseteq G$ invariably generates $G$
if and only if $S \not\subseteq \widetilde{H}$ for all
proper subgroups $H < G$. If $G$ is finitely generated
then $S \subseteq G$ invariably generates $G$ if and
only if $S \not\subseteq \widetilde{M}$ for all maximal
subgroups $M < G$.
\end{lemma}

This implies the following easy observation.

\begin{lemma} The following are equivalent for a group $G$.

(i) $G$ is IG.

(ii) For every proper subgroup $H < G$ we have $\widetilde{H} \neq G$.

(iii) If $H \le G$ and $H$ intersects every conjugacy class
of $G$ then $H=G$.

(iv) In every transitive action of $G$ on a set $X$ with more than
one element there is $g \in G$ acting on $X$ as a fixed-point-free
permutation.

\end{lemma}

It is also easy to see, more generally, that $S \subseteq G$ generates $G$
invariably if and only if in any transitive action of $G$
on a set with more than one element there exists $s \in S$
acting fixed-point-freely.

Using Lemma 2.2 we readily see that finite groups are IG.
Groups satisfying condition (iv) above were studied by
Wiegold and others, see \cite{W1, W2, CSW}. Reformulating results
from \cite{W1, W2} using Lemma 2.2 we obtain the following.

\begin{Corollary}

(i) Virtually solvable groups are IG.

(ii) Nonabelian free groups are not IG.

(iii) The class of IG groups is extension closed.

(iv) The class of IG groups is not subgroup closed.

\end{Corollary}

Other examples of non IG groups are infinite groups $G$ all of
whose nontrivial elements are conjugate (see \cite{HNN} and
\cite{Os}); indeed in such groups we have $\widetilde{H}=G$
for every nontrivial subgroup $H < G$.

A wide class of algebraic groups is also not IG. Indeed
we have the following characterization.

\begin{Proposition} Let $G$ be a linear algebraic group
over an algebraically closed field. Then $G$ is IG
if and only if it is virtually solvable.
\end{Proposition}

\proof
If $G$ is virtually solvable then it is IG by Corollary~2.3.

Now suppose $G$ is IG.
By a theorem of Steinberg (see Theorem 7.2 of \cite{St}),
every automorphism of a linear algebraic group $G$ fixes
some Borel subgroup of $G$. This implies that
if $g$ is any element of $G$, then $B^g = B$ for some Borel subgroup
$B$ of $G$. Thus the union of the normalizers $N_G(B)$ over the Borel
subgroups $B$ of $G$ equals $G$. Since the Borel subgroups are all conjugate,
it follows that
$$
\widetilde{N_G(B)} = G
$$
for any Borel subgroup $B$ of $G$. Lemma 2.1 and the assumption
that $G$ is IG now imply that $N_G(B) = G$. This in turn implies
that $G$ is virtually solvable.
\qed
\medskip

Let $\Phi(G)$ denote the Frattini subgroup of a finitely generated
group $G$. Then a subset of $G$ generates $G$ if and only if its image
in $G/\Phi(G)$ generates $G/\Phi(G)$.
It follows that a subset of $G$
invariably generates $G$ if and only if its image in
$G/\Phi(G)$ invariably generates $G/\Phi(G)$.


For an FIG group $G$, let $d_I(G)$ denote the minimal
number of invariable generators for $G$.

\begin{Lemma}
Let $G$ be a finitely generated group.

(i) If $G/\Phi(G)$ is IG, then so is $G$.

(ii) If $G/\Phi(G)$ is FIG, then so is $G$.

(iii) If $G/\Phi(G)$ is finite, then $G$ is FIG.

(iv) $d_I(G) = d_I(G/\Phi(G))$.

(v) If $G/\Phi(G)$ is a finite (nonabelian) simple group, then $d_I(G)=2$.

\end{Lemma}

\proof Parts (i)-(iv) follow immediately from the remarks
preceding the lemma.
Part (v) follows from (iv) and a result from \cite{KLS}:  finite simple groups are invariably generated by two
elements.
\qed
\medskip


\begin{Lemma} Let $G$ be a finitely generated group.

(i) If all maximal subgroups of $G$ have finite index
then $G$ is IG.

(ii) Suppose that there exists an integer $c$ such that every maximal
subgroup $M$ of $G$ satisfies $|G:M| \le c$. Then $G$ is FIG.

\end{Lemma}

\proof If $M < G$ has finite index then $\widetilde{M} \ne G$.
Part (i) now follows from Lemma 2.1.

To prove part (ii), note that it follows from the assumption on $G$
that $G$ has finitely many maximal subgroups, and since they all have
bounded index we see that $G/\Phi(G)$ is finite. The result now follows
from part (iii) of Lemma 2.5.
\qed
\medskip


We now apply the lemma to the Grigorchuk group   \cite{Gr}:

\begin{Corollary}
The Grigorchuk group $G$ is FIG. In fact $d_I(G) = 3$.
Thus, residually finite FIG groups need not be virtually solvable.
\end{Corollary}

\proof Recall that the Grigorchuk group $G$ is an infinite $2$-group
generated by $3$ elements (of order $2$).
It was shown by Pervova in \cite{Pe} that all maximal subgroups of
this group have finite index, hence they are of index $2$.
 The result follows from Lemma 2.6. In fact $G/\Phi(G)$ is an
elementary abelian group of order $8$ and hence, by 2.5(iv),
$G$ is invariably generated by $3$ elements.
\qed
\medskip

We continue with additional basic  results on FIG groups.


\begin{lemma} The class of FIG groups is extension closed.
\end{lemma}

\proof Let $N \lhd G$ and suppose both $N$ and $G/N$ are
FIG. We need to show that $G$ is FIG. Suppose $N$ is invariably
generated by a finite set $S$, and $G/N$ is invariably generated
by a finite set $T$. Let $T_1 \subseteq G$ be a set of representatives
for $T$ in $G$.

We claim that $G$ is invariably generated by the finite set $S \cup T_1$.
To show this, let $H \le G$ with $\widetilde{H} \supseteq S \cup T_1$.
We need to show that $H = G$. Clearly $\widetilde{HN/N} \supseteq T$,
which implies $HN/N =G/N$ so $HN = G$.

Now let $s \in S$. Then there exist $h \in H$ and $g \in G$ such that
$s = h^g$. Write $g = h_1 n$ where $h_1 \in H$ and $n \in N$.
Then $s = h^{h_1n} = h_2^n$ where $h_2 \in H$. Since $s \in N$
it follows that $h_2 \in H \cap N$. Thus $S$ is covered by the
union of $N$-conjugates of $H \cap N$. This implies
$H \cap N = N$, so $H \supseteq N$. Therefore $H = G$.

\qed
\medskip


\begin{Corollary} Suppose $N \lhd G$ has finite index, and $N$ is
FIG. Then $G$ is FIG.
\end{Corollary}

\proof This follows from Lemma 2.8 above, since $G/N$ is finite,
hence FIG.

\qed
\medskip


\begin{lemma} Let $A \lhd G$ be an abelian normal subgroup.
Suppose $G/A$ is FIG, and $A$ is finitely generated as a
$G/A$-module. Then $G$ is FIG.
\end{lemma}

\proof Let $S \subseteq A$ be a finite set generating $A$
as a $G/A$-module. Let $T \subseteq G/A$ be a finite set
which invariably generates $G/A$, and let $T_1$ be a
a set of representatives for $T$ in $G$. We claim that
the finite set $S \cup T_1$ invariably generates $G$.

Indeed, let $H \le G$ with $\widetilde{H} \supseteq S \cup T_1$.
We have to show that $H = G$.
As in the proof of Lemma 2.8, we obtain $HA = G$.
Let $s \in S$ and $g \in G$. Then $s = h_0^{g_0}$ for some $g_0 \in G$
and $g_0g = h_1a_1$ for some $h_1 \in H$ and $a_1 \in A$.
Hence $s^g = h_0^{g_0g} = h_0^{h_1a_1} = h_2^{a_1}$ for some
$h_2 \in H$. But $s^g \in A$, so $h_2 \in H \cap A$ and since $A$
is abelian we have $s^g = h_2^{a_1^{-1}} = h_2  \in H \cap A$.

However, the elements $s^g$ ($s \in S, g \in G$) generate $A$.
It follows that $H \cap A = A$, so $H \supseteq A$ and $G = H$
as required.

\qed
\medskip

We can now derive some consequences.


\begin{Proposition} Let $G$ be a finitely generated group.

(i) If $G$ is a solvable Max-n group then $G$ is FIG.

(ii) If $G$ is polycyclic then $G$ is FIG.

(iii) If $G$ is abelian-by-polycyclic then $G$ is FIG.

(iv) If $G$ is abelian-by-nilpotent then $G$ is FIG.

(v) If $A_1, \ldots , A_m$ are finitely generated abelian
groups, then the iterated wreath product
$A_1 \wr (A_2 \wr (A_3 \wr \ldots \wr A_m)))$
is FIG. In particular, the lamplighter group $C_2 \wr \Z$
is FIG.

\end{Proposition}

\proof
Recall that a group $G$ is a Max-n group if it satisfies the maximal
condition on normal subgroups. This is equivalent to every
normal subgroup of $G$ being finitely generated as a normal subgroup.

We   prove part (i)   by
induction on the derived length $d$ of $G$. The result is clear
for abelian groups, so suppose $d > 1$. Let $A = G^{(d-1)}$.
Then $A \lhd G$ is abelian and finitely generated as a normal
subgroup. By induction hypothesis, $G/A$ is FIG, so  $G$ is FIG
by Lemma 2.10.

It is well known (see \cite{H}) that polycyclic groups
and finitely generated abelian-by-polycyclic groups -- and in particular
abelian-by-nilpotent groups -- are Max-n.
Thus parts (ii)-(iv) follow.

Part (v) is proved by induction on $m$ using Lemma 2.10.
\qed
\medskip

Of course if $G$ has a finite index subgroup satisfying
one of the above conditions (i)-(v) then it is also FIG.

It is known that a finitely generated center-by-metabelian
group need not be Max-n, indeed its center need not be
finitely generated \cite{H}.


\begin{Proposition}

(i) Let $G$ be a group and $N \lhd G$ a nilpotent normal subgroup.
Suppose $G/N$ is FIG and $N$ is finitely generated as a normal subgroup.
Then $G$ is FIG.

(ii) A finitely generated metanilpotent-by-finite group is FIG.
\end{Proposition}

\proof

We need the following.
\medskip

{\bf Claim.} Let $G$ be a group, $N \lhd G$ a nilpotent normal subgroup.
If $G/N'$ is FIG then so is $G$.

To show this suppose $G/N'$ is invariably generated by the finite subset $S$,
and let $S_1 \subseteq G$ be a set of representatives for $S$ in $G$.
We claim that $S_1$ invariably generates $G$. To show this, let
$H \le G$ such that $\widetilde{H} \supseteq S_1$ and conclude that
$H = G$.

Since $HN'/N' \supseteq S$ we have $HN' = G$. If $n \in N$ then
$n = hn'$ for some $h \in H$ and $n' \in N'$. Thus $h = n n'^{-1}$
so $h  \in H \cap N$. It follows that $(H \cap N)N' = N$.

It is well know that if $L$ is a subgroup of a nilpotent group $N$ satisfying
$L N' = N$ then $L = N$. Applying this for $L = H \cap N$ we obtain
$H \cap N = N$, so $H \supseteq N$. But $HN = G$, hence $H = G$,
proving the claim.
\medskip

Next, we prove part (i). By Lemma 2.10, $G/N'$ is FIG.
Hence, by the claim above, $G$ is FIG.

To prove part (ii) apply Corollary 2.9 to reduce to the case
when $G$ is metanilpotent.
Let $N \lhd G$ such that $N$ and $G/N$ are nilpotent.
Then $G/N'$ is abelian-by-nilpotent, hence it is FIG by Proposition 2.11(iv).
It now follows from the claim above that $G$ is FIG.

\qed
\medskip

We see from Proposition 2.12 and the remark preceding it that
finitely generated solvable groups which are FIG need not
satisfy Max-n.
It is also easy to see using the arguments above that an iterated
wreath product of finitely generated nilpotent groups is FIG.

\section{Finitely generated linear groups}

In this section we prove Theorem 1.1.

Let $G$ be a linear group. If $G$ is finitely generated and virtually
solvable, then, by the Lie Kolchin Theorem, $G$ contains a finite index
subgroup represented (up to conjugacy) by upper triangular matrices.
Hence $G$ is nilpotent-by-abelian-by-finite, so it is FIG by
Proposition 2.12.

To prove the other direction we will assume $G$ is a subgroup of $GL_n(F)$
for some field $F$, and that it is FIG but not virtually solvable.
We will derive a contradiction by using the Strong Approximation Theorem
(see \cite{We, P, N} and page 406 in \cite{LS}).

Since $G$ is finitely generated, it is   contained in $GL_n(A)$
for some
finitely generated subring $A$ of $F$. By Theorem 4.1 of \cite{LL}
there exists
a specialization, namely a ring homomorphism $\phi: A \rightarrow k$, where
$k$ is a global field, such that the image of $G$ under the induced
map $\phi_1: GL_n(A) \rightarrow GL_n(k)$ is not virtually solvable. Replacing
$G$ by $\phi_1(G)$ we shall assume $F = k$ (and $G$ is still FIG
as a quotient of an FIG group).

Let $H$ be the Zariski closure of $G$ in $GL_n(\overline{k})$,
where $\overline{k}$ is the algebraic closure of $k$.
Then $H$ is a linear algebraic group (over an algebraically
closed field) which is not virtually solvable. Dividing $H$ by
its maximal solvable normal subgroup we can assume that $H$
is semisimple. Furthermore,
by factoring out a suitable normal subgroup we may assume
that $H$ is homogeneous of the form $L^m \rtimes \Delta$ where
$L$ is a simple algebraic group of adjoint type and $\Delta$ is a finite group
(permuting the copies of $L$ transitively and possibly acting as
outer automorphisms  on each copy).
The image of $G$ in this process is still FIG, not virtually solvable,
and Zariski dense. We replace $G$ by this image.

Let $L_1$ be the simply connected cover of $L$ and let
$\psi: L_1 \rightarrow L$ be the covering map. The finite group $\Delta$
acts also on $L_1^m$ and we obtain an epimorphism
$\psi_1:H_1: = L_1^m \rtimes \Delta
\rightarrow H = L^m \rtimes \Delta$.
The group $\psi_1^{-1}(G)$ is a central extension of $G$ with a
finite center, and hence is also $FIG$ by Lemma 2.8. Replacing $H$ by
$H_1$ and $G$ by $G_1$ we can assume that $G$ is an FIG dense subgroup
of an algebraic group $H \le GL_{n_1}$ whose connected component $H^0$
is simply connected. Furthermore, by restriction of scalars we can even assume
that $G$ is inside $GL_n(k)$ for some $n$, where $k = \Q$ or $\F_p(t)$.
Moreover, $G$ is inside $H(O_T)$, where $O$ is the ring of integers
of $k$  and $T$ is a finite set of primes.

We are now in a position to apply the Strong Approximation Theorem.
According to this theorem there exists a finitely generated
ring $R$ of $O_T$ such that $k$ is the field of fractions of $R$ (in
characteristic $p$ this may require replacing the original field $k$
by a smaller subfield), $G$ is inside $H(R)$ and, for almost
every prime ideal $P$ of $R$,  the image of $G^0 = G \cap H^0$ in
$H^0(R/P)$ is onto. Note that for almost every prime $P$, $H^0(R/P)$
and $H(R/P)$ are
well defined, as $k$ is the ring of fractions of $R$, and $H, H^0$
are both defined over $k$, since $G \le H(R)$ is Zariski dense in $H$.
Moreover, $H^0(R/P) \rtimes \Delta$ is also well defined.
Since
$G^0$ is mapped onto $H(R/P)$ and $G$ is mapped onto $H/H_0$, $G$ is mapped onto $H^0(R/P) \rtimes \Delta$.

Now let $S \subset G$ be a finite set which invariably generates $G$.
By Proposition 2.4 and its proof, for each $s \in S$ there exists
an element $h(s) \in H$ such that $s^{h(s)} \in B_1 = N_H(B)$,
where $B$ is a Borel subgroup of $H^0$. Note that $B_1$ is virtually
solvable. The finitely many elements $h(s), s \in S$ all belong to
$H(k_1)$, where $k_1$ is a finite extension of $k$. In fact they
are even in $H(O^1_{T_1})$ where $O^1$ is the ring of integers of
$k_1$ and $T_1$ is a finite set of primes. By extending $T_1$ further
if needed we may assume $R_1: = O^1_{T_1} \supseteq R$.

By the Chebotarev density theorem there exist infinitely many prime
ideals $P$ of $R$   that split completely in $k_1$; in particular
for such $P$, there exists a prime ideal $P_1$ of $R_1$ for which
the inclusion $R \subseteq R_1$ induces an isomorphism
$R/P \cong R_1/P_1$. For almost all such primes $P$ the image
of $G$ in $H(R_1/P_1) \cong H(R/P)$ is onto, while the image $B_2$
of $G \cap B_1$ there is a proper subgroup, since this image
is solvable-by-bounded. The image of $h(s)$ there conjugates $s$
into $B_2$. Therefore $H(R/P)$, the finite quotient of $G$,
is not invariably generated by the image of $S$.
This contradiction completes the proof of Theorem 1.1.
\qed
\medskip

{\bf Remark.} Our proof in fact shows something stronger: for every
non-virtually solvable linear group $G$ and every finite subset
$S$ of it, there exists a proper {\em finite index} subgroup $H < G$
such that $\widetilde{H} \supseteq S$. Clearly,
for such $H$, $\widetilde{H} \ne G$. It is possible that there
exists an infinite index subgroup $H$ with  $\widetilde{H} = G$.
For example, this happens in (nonabelian) free groups $G$. But we do not know
if this is the case for all non virtually solvable linear groups,
i.e., whether there exists a linear IG group which is not virtually
solvable.

\section{Profinite groups}

Let $G$ be a profinite group. Then generation and invariable generation
in $G$ are interpreted topologically, and by subgroups we mean closed
subgroups. It is then easy to see that the basic results in Section 2
also hold in the category of profinite groups.

Just as every finite group is IG, every profinite group $G$ is also
IG. Indeed every proper subgroup of a profinite group $G$ is
contained in a maximal open subgroup $M$, and since $M$ has finite
index we have $\widetilde{M} \ne G$. Hence $G$ is IG by Lemma 1.1.

On the other hand, finitely generated profinite groups are not
necessarily FIG. In fact in Proposition 2.5 of \cite{KLS} we showed
that there exist $2$-generated finite groups $H$ with $d_I(H)$
(the minimal number of invariable generators)
arbitrarily large. This implies that the free profinite
group ${\widehat{F_d}}$ on $d \ge 2$ generators is not FIG.

On the other hand, the free pro-$p$ group on $d < \infty$ generators
is FIG, since its Frattini subgroup is of finite index (see Lemma 2.5(iii)
above). Since free pronilpotent groups are direct products of free
pro-$p$ groups, we easily deduce that every finitely generated
pronilpotent group is FIG. Compare this with Problem 4 in Section 5 below
regarding prosolvable groups.

The following lemma is useful in the proofs of Proposition~1.3.


\begin{Lemma} Let $G$ be a simply connected simple Chevalley group.

(i) $\Phi(G(\Z_p))$ contains the second congruence subgroup.

(ii) If $p > 3$ then $\Phi(G(\Z_p))$ is the first congruence subgroup.

(iii) For a prime power $q$ (with finitely many possible exceptions),
$\Phi(G(F_q[[t]])$ is the second congruence subgroup.

\end{Lemma}

\proof See \cite{Wei} and \cite{LL}.
\qed
\medskip

{\em Proof of Proposition 1.3.}  Recall that $G$ is a simply connected
simple Chevalley group. It is well known that the profinite group
$G(\Z_p)$ has an open finitely generated pro-$p$ subgroup.
This implies that its Frattini
subgroup $N$ is open. Using part (iii) of 2.5 we see that $G(\Z_p)$
is FIG. Moreover, by Lemma 4.1(i) we see that the Frattini quotient $Q$
of $G(\Z_p)$ is a finite simple group, or
an extension of an abelian group $A$ by a finite
(quasi-)simple group $T$. Moreover, in the latter case, $A$ is
generated as a normal subgroup by a single element.

Thus, in the first case we have $d_I(G(\Z_p)) = d_I(Q) = 2$
by Lemma~2.5(v), while in the second case we have $d_I(G) = d_I(Q) \le 3$
by Lemma 2.8 and its proof.

Hence in any case $G(\Z_p)$ is invariably generated by $3$ elements,
which we denote by $g_1(p), g_2(p) , g_3(p)$.

Now, the adelic group $G(\widehat{\Z})$ is isomorphic to the direct product
$\prod_p G(\Z_p)$. For $i = 1, 2, 3$ let $g_i$ denote the sequence
$(g_i(p))$ where $p$ ranges over the primes. Then it is easy to see
that $g_1, g_2, g_3$ generate $G(\widehat{\Z})$ invariably.
This proves part (i) of the Proposition.

Next, if $p>3$, then by part (ii) of Lemma 4.1, the Frattini quotient
of $G(\Z_p)$ is a finite simple group, so (as argued above)
$G(\Z_p)$ is invariably generated by two elements. This proves part (ii).
\qed
\medskip

We next  generalize Proposition 1.3 and deal with groups
over arbitrary global fields. This requires some preparations.


\begin{Lemma}
Let $G=T^m$ for a nonabelian finite simple group $T$.
Let $S=\{s_1,\dots,s_r\}\subset G$, so that
$s_i=(t_1^i,\dots,t_m^i), t_j^i\in T$.  Form the matrix
$$
A=\begin{pmatrix} 
        t_1^1&\dots & t_m^1 \\
         &\dots & \\
        t_1^r&\dots & t_m^r \\
     \end{pmatrix}.
     $$
Then $S$ invariably generates $G$ if and only if the following both hold$:$
\begin{itemize}
\item[\rm(a)]
If $1\le j\le m$ then $\{ t_j^1,\dots,t_j^r  \}$ generates $T$ invariably.
\item[\rm(b)] The columns of $A$ are in
different $\Aut(T)$-orbits for the diagonal action of  $\Aut(T)$ on $T^r$.
\end{itemize}
\end{Lemma}

\proof This follows immediately from the generation criterion for $T^m$
in \cite[Proposition 6]{KL}.
\qed
\medskip

The number of conjugacy classes of a finite group $T$ is denoted
by $k(T)$. The next result shows that rather large powers of finite
simple groups are still invariably generated by few elements.


\begin{Proposition} Let $T$ be a finite simple group.
Given $r \ge 2$, let $m(T,r)$ denote the maximal integer $m$
such that $d_I(T^m) \le r$. Then
$$
k(T)^{r-2}/|Out(T)|-1 < m(T,r) \le k(T)^r.
$$
\end{Proposition}

\proof
 Suppose $T$ is invariably generated by $a, b \in T$.
Let $A, B \subset T$ be the conjugacy classes of $a, b$ respectively.
Consider all  $r$-tuples $(A, B, C_3, \ldots , C_r)$ where each $C_i$
ranges over all conjugacy classes of $T$. There are $k(T)^{r-2}$ such
tuples, and they split into at least $x: = k(T)^{r-2}/|Out(T)|$
different orbits under the action of $Out(T)$. Therefore
if $m$ is the greatest integer in $x$ then it follows from Lemma~4.2
that $T^m$ is invariably generated by $r$ elements, two
of which are $(a, a, \ldots , a)$, $(b, b, \ldots , b)$.
This proves the lower bound on $m(T,r)$.

The upper bound follows immediately from Lemma~4.2.
\qed
\medskip

Note that $|Out(T)| \le \log |T|$,  whereas $k(T)$ is much larger:
  it is roughly $c^{\sqrt{n}}$ if $T = A_n$ and $q^l$ if $T = G(q)$,
a Lie type group of rank $l$ over the field with $q$ elements (see
\cite{FG}).
This shows that the lower and upper bounds in Proposition 4.3 are
of rather similar orders of magnitude.


\begin{Corollary}

(i) If $m \in \N$ satisfies $m \le k(T)/|Out(T)|$ then $d_I(T^m) \le 3$.

(ii) For every $m \in \N$ and almost all finite simple groups $T$
we have $d_I(T^m) \le 3$.

(iii) Let $G$ be a Chevalley group and $c \in \N$ a given constant.
Then for all sufficiently large prime powers $q$ we have
$d_I(G(q)^{cq}) \le 4$.

(iv) Let $a_n \in \N$ be such that $\log a_n / \sqrt{n} \rightarrow \infty$
as $n \rightarrow \infty$. Then $d_I(A_n^{a_n}) \rightarrow \infty$
as $n \rightarrow \infty$.

\end{Corollary}

\proof
Part (i) follows immediately from 4.3.

Part (ii) follows from (i) and the remark above, implying that
$k(T)/|Out(T)| \rightarrow \infty$ as $T$ ranges over the
finite simple groups.

For part (iii), we easily verify using \cite{FG} that
$k(G(q))^2/|Out(G(q))| \ge cq$ if $q$ is sufficiently large
(given $c$). Using Proposition~4.3 with $r = 4$ yields the result.

Part (iv) follows from the upper bound in Proposition~4.3.
\qed
\medskip

We can now prove the main result leading to Theorem 1.4.


\begin{theorem}
Let $k$ be a global field and $T$ a finite set of places of $k$
containing all the archimedean ones. Let $G$ be a connected
simply connected simple algebraic $k$-subgroup of $GL_n$. Let
$\A_T = \prod^*_{v \not\in T} k_v$ be the ring of $T$-adeles of $k$,
and let $H$ be an open compact subgroup of $G(\A_T)$.
Then $H$ is an FIG profinite group.
\end{theorem}

\proof
 The structure of the proof is similar to that of Proposition 1.3,
but there are more technicalities to handle. As shown in the proof
of Theorem 3.1 in \cite{LL}, after passing to a finite index
subgroup, $H$ is the product of infinitely many groups $H_v$, where
$H_v$ is a virtually pro-$p$ open subgroup of $G(k_v)$ for the
various completions $k_v$, $v \not\in T,$ of $k$.

Factoring out the Frattini subgroup of $H$, we are left with
an infinite product of finite groups.
For almost every $v$,
$H_v/\Phi(H_v) $ is an extension of a finite elementary abelian group $M_v$
generated as a normal subgroup by boundedly many elements by a finite
(quasi)simple group $T_v$ of the same type as $G$ over $F_v := O_v/m_v$,
the residue field of $k_v$.
In fact, with the exception of finitely many
group types and characteristics, $M_v$ is abelian and simple as a
$T_v$-module, hence generated as a normal subgroup by one element;
moreover if $char(k)=0$ then $M_v = 0$.
See the proof of [LL, Theorem~3.1] and especially properties
(a), (b), (c) there.
We may ignore finitely many factors.

Now, a simple group of Lie type $G$ over a finite field of order
$q$ occurs in this product with bounded multiplicity if
$char(k) = 0$ and with multiplicity $\le cq$ (for some constant
$c$) if $char(k) > 0$.
So, in any case, $T:= \prod T_v$ is FIG
by Corollary~4.4(iii). Moreover, $M = \prod M_v$ is generated
by boundedly many elements as a normal subgroup.
Hence, by part (i) of Proposition 2.12, $H/\Phi(H)$ is FIG, and
so is $H$ by 2.5(ii).

\qed
\medskip

{\it Proof of Theorem 1.4. }
Since $G(O_T)$ has CSP,
its profinite completion is an extension of a finite center
by a group $H$ as in Theorem 4.5. The result follows
from Theorem~4.5 and Lemma~2.8.
\qed
\medskip

We now make preparations for the proof of Theorem 1.5.
For background on Fuchsian groups, see \cite{LiSh} and the
references therein.

Higman conjectured that if $G$ is any Fuchsian group, then every large
enough alternating group $A_n$ is a quotient of $G$. This was proved
in \cite{E} (in the oriented case) and \cite{LiSh} provides a
probabilistic proof
of the conjecture (also in the non-oriented case).
In fact the following strengthening of Higman's
conjecture also holds.


\begin{Proposition} Let $G$ be a Fuchsian group (oriented or non-oriented).
If $n$ is sufficiently large, and $a_n$
is the integral part of $(n!)^{1/43}$, then $A_n^{a_n}$ is
a quotient of $G$.
\end{Proposition}

\proof
Let $\mu(G)$ denote the measure of $G$, namely $- \chi(G)$, where
$\chi(G)$ is the Euler characteristic of $G$.
It is known that $\mu(G) \ge 1/42$. By Theorem 1.1 of \cite{LiSh}
and the remark following it we have
$$
|Hom(G, A_n)| \ge (n!)^{\mu(G) + 1 + o(1)} \ge (n!)^{43/42 + o(1)}.
$$
By Theorem 1.7 of \cite{LiSh}, most of the homomorphisms from
$G$ to $A_n$ are epimorphisms, and so
$$
|Epi(G, A_n)| \ge (n!)^{43/42 + o(1)},
$$
where $Epi(G, A_n)$ is the set of epimorphisms from $G$ to $A_n$.

Suppose $G$ is generated by $g_1, \ldots , g_r$. Every epimorphism
$\phi: G \rightarrow A_n$ gives rise to an $r$-tuple
$(\phi(g_1), \ldots , \phi(g_r)) \in A_n^r$ which generates $A_n$.
Form a matrix whose columns are these $r$-tuples for all
$\phi \in Epi(G, A_n)$. Let $S_n = Aut(A_n)$ act on these $r$-tuples
diagonally. Then there are at least
$|Epi(G, A_n)|/|S_n| \ge (n!)^{1/42+o(1)}$ different orbits under this
action. Since $a_n \le (n!)^{1/43}$ it follows using
\cite[Proposition 6]{KL} that $A_n^{a_n}$ is a quotient of $G$.
\qed
\medskip


\begin{Lemma} Let $a_n$ be as above. Then
$d_I(A_n^{a_n}) \rightarrow \infty$ as $n \rightarrow \infty$.
\end{Lemma}

\proof This follows from part (iv) of Corollary 4.4.
\qed
\medskip

{\it Proof of Theorem 1.5.}
The theorem follows immediately from
Proposition~4.6 and Lemma~4.7.
\qed
\medskip

\section{Open problems}

We conclude this paper by posing some natural problems
which may inspire further research.
\medskip

1. Is a finite index subgroup of an IG group necessarily IG?
\medskip

2. Is a finite index subgroup of an FIG group necessarily FIG?
\medskip

3. Are finitely generated solvable groups FIG?
\medskip

4. Are finitely generated prosolvable groups FIG?
\medskip

5. Are finitely generated solvable profinite groups FIG?
\medskip

6. Is $SL_n(\Z)$ ($n \ge 3$) IG?
\medskip

7. Is $SL_n(\Q)$ IG?
\medskip

8. Is every IG linear group virtually solvable?
\medskip

9. Is every (non-elementary) word hyperbolic group non IG?
\medskip

10. Is the profinite completion of every (non-elementary) word hyperbolic
group non FIG?

\end{document}